\documentclass{amsart}
\usepackage{fullpage}
\usepackage{amsmath,amssymb,amsthm}
\usepackage{enumerate}
\usepackage{subfigure}
\usepackage{tikz}
\usetikzlibrary{intersections}
\usetikzlibrary{3d}
\usetikzlibrary{calc}
\usetikzlibrary{decorations.pathreplacing}
\tikzstyle{dot}=[circle, fill=black,inner sep=0.5mm]
\tikzstyle{line}=[thick]

\begin{document}
\newtheorem{lemma}{Lemma}
\newtheorem{theorem}{Theorem}
\newtheorem{question}{Question}
\newtheorem{conjecture}{Conjecture}

\newcommand{\trip}[1]{[#1]}
\newcommand{\sv}[1]{y^{#1}}
\newcommand{\tspace}[1]{X^{#1}}
\newcommand{\ignore}[1]{}

\newcommand{\affhull}{\text{aff.hull}}
\newcommand{\linhull}{\text{lin.hull}}
\newcommand{\QQ}{\mathcal{Q}}
\newcommand{\UU}{\mathcal{U}}
\newcommand{\CC}{\mathcal{C}}
\newcommand{\cl}{\mbox{cl}}
\newcommand{\R}{\mathbb{R}}
\newcommand{\F}{\mathbb{F}}
\newcommand{\GF}{{\rm GF}}
\newcommand{\supp}{\text{supp}}

\title{Reconstructing a phylogenetic level-1 network from quartets}
\begin{abstract} We describe a method that will reconstruct an unrooted binary phylogenetic  level-1 network on $n$ taxa from the set of all quartets containing a certain fixed taxon, in $O(n^3)$ time.
We also present a more general method which can handle more diverse quartet data, but which takes $O(n^6)$ time.
Both methods proceed by solving a certain system of linear equations over $\GF(2)$.

For a general dense quartet set (containing at least one quartet on every four taxa) our $O(n^6)$ algorithm constructs
a phylogenetic level-1 network consistent with the quartet set if such a network exists and returns an ($O(n^2)$ sized) 
certificate of inconsistency otherwise. This answers a question raised by Gambette, Berry and Paul regarding the
complexity of reconstructing a level-1 network  from a dense quartet set.
\end{abstract}

\author{J.C.M. Keijsper and R.A. Pendavingh}

\maketitle

\section{\label{sect:intro}Introduction}
In phylogenetics, a rooted or  unrooted phylogenetic tree is traditionally  used as a model for
studying evolutionary relationships. In addition, there has been a recent interest in developing methods for modelling and reconstructing reticulation. A rooted or unrooted phylogenetic network can be an appropriate model in case of reticulate evolution. 
 In a rooted context, a level parameter was introduced to measure the level of complexity  of a rooted phylogenetic network
\cite{JanssonSung2006}, and polynomial time algorithms were proposed to reconstruct rooted phylogenetic networks of bounded level from a
set of triplets (rooted subtrees on three taxa) \cite{JanssonNguyenSung2006}, \cite{IerselKeijsperKelkStougie2009}, \cite{ToHabib2009}, \cite{IerselKelk2011} in case the triplet set is dense (i.e. contains a triplet on every subset of three taxa).
In \cite{Gambette2012}, an unrooted analogue of this level parameter was introduced, but it was left as an open problem whether an unrooted level-1 network can be reconstructed in a similar manner from a dense set of unrooted subtrees on four taxa in polynomial time. Note that
for general (nondense) quartet sets, the decision problem whether a compatible tree  or level-1 network exists is NP-complete \cite{Steel1992}, \cite{Gambette2012}.

This paper deals with the problem of recovering an unrooted binary phylogenetic level-1 network $G$ 
from  a set of given {\em quartets}, each quartet describing a subtree of $G$ on a four-tuple of taxa.
It is known that  any binary tree $G$ whose leaves are labelled by $X$, where $|X|=n$, can be recovered in polynomial ($O(n^4)$) time from the full set of quartets determined by $G$, see e.g. \cite{BerryGascuel2000}, or \cite{SempleSteel2003}. Gambette, Berry and Paul \cite{Gambette2012} consider the related problem of reconstructing a level-1 network from quartets, and they give an $O(n^4)$ algorithm that  reconstructs a level-1 network $G$ from the full set of quartets induced by $G$.

In this paper, we describe an algorithm which, given the subset of quartets involving some fixed taxon $\infty\in X$, will recover a level-1-network $G$ exhibiting exactly these quartets on $\infty$, if such a network exists, in $O(n^3)$ time.
 The returned level-1 network is as `sparse', or `tree-like' as possible, in the sense that it has the maximum number of 
 cut-edges among all 
 compatible level-1 networks.
 We also present a slower ($O(n^6)$) but more flexible method, which may take any set of quartets $Q$, and which terminates with either a sparsest compatible level-1 network, or a subset $Q'\subset Q$ of $O(n^2)$ quartets which is already incompatible with any level-1 network, or advice where to add further quartets, if the given set of quartets was insufficient for the method to reach a definite conclusion.
 A dense set of quartets is always sufficient for the latter method to conclude incompatibility or find a sparsest compatible level-1 network, which implies that a polynomial algorithm exists for reconstructing a level-1 network from a dense set of quartets.
 We sketch these methods here, noting that the more standard definitions on phylogenetic networks can be found in the next section.

Unrooted level-1 networks are {\em outerplanar}, i.e. can be drawn in a disk without crossings such that all the vertices appear on the boundary of the disk. Any such drawing determines a cyclic ordering of the set of taxa at the leaves of such a network. We will say that  a cyclic ordering of $X$ is {\em compatible} with a network $G$ on $X$ if it arises from  an outerplanar drawing of $G$ in this manner.  A level-1 network $G$ on $X$ may be compatible with many cyclic orderings of $X$. Consider the fact that if $e$ is any cut-edge of $G$, then an outerplanar drawing of $G$ can be modified  by `flipping' the drawing of one component of $G\setminus e$. If one of the components of $G\setminus e$ is labelled by the set of taxa $Y\subseteq X$, then $Y$ and $X\setminus Y$ will appear consecutively in any cyclic ordering compatible with $G$, and flipping the drawing as described will result in a cyclic ordering in which the ordering of one of $Y$ or $X\setminus Y$ is reversed.
Combinatorial properties of cyclic orderings compatible with trees were studied before in \cite{SempleSteel2004}.

We will argue below that if $G$ is a level-1 network on $X$, then $G$ is fully determined by the complete set of cyclic orderings compatible with $G$.
The set of cyclic orderings itself can be recovered from the full set of quartets compatible with $G$ --- in fact, we will argue that the set of quartets containing some fixed taxon $\infty\in X$  suffices. 

Calculating orderings from the quartets, and then the network from the orderings, we obtain an $O(n^3)$ algorithm for recovering $G$ from the set of  all quartets containing a fixed taxon $\infty\in X$.  A key observation in this paper is that the set of all compatible cyclic orderings, upon suitably encoding each cyclic ordering as a $0,1$-vector, determines an affine space  over $\GF(2)$. This affine space is exactly the set of solutions of a system of linear equations, with one equation for each given quartet. In its simplest form, our algorithm proceeds by drawing up these $O(n^3)$ linear equations (over some $O(n^2)$ variables) and solving them. From one solution of this system, which encodes  some cyclic ordering compatible with $G$, we then go on to recover $G$.
Interestingly, the complete set of solutions to this system of equations (given by a particular solution and a basis of the parallel linear subspace) can be viewed as a compact encoding of the set of all possible level-1 networks compatible with the given quartet set. 
Thus, our method does not merely return one possible solution, but in fact it can be made to return a (polynomially sized) representation of all possible solutions.

Our system of $O(n^3)$ linear equations in $O(n^2)$ variables is sparse, in the sense that each equation involves at most 4 variables, and with a suitable choice of variables relative to $\infty$, just two variables. In the latter case, finding a solution to the system is about as hard as finding a spanning forest in a graph with $O(n^2)$ vertices and $O(n^3)$ edges, and takes $O(n^3)$ time.

Without the assumption that all given quartets involve $\infty$, we lose the special structure of the system of equations. But we will still have a sparse system of equations with at most 4 nonzero variables in each equation. Using standard  linear algebra for solving these equations, we obtain a 
general method for inferring level-1 networks from quartets. 
Given any inconsistent set of quartets, the algorithm will be able to produce a subset of at most $O(n^2)$ of the input quartets which is already inconsistent with any level-1 network. If the presented quartets do not suffice for the method to determine $G$, the method returns 4-tuples of taxa so that any new quartet data concerning one of these 4-tuples will make it possible for the method to advance. Used interactively in this manner, the method will determine a level-1 network in $O(n^2)$ steps, taking $O(n^4)$ time in each step.
So the total running time of this more flexible method is worse at $O(n^6)$ time, due to the use of more general solution methods for solving linear equations over $\GF(2)$, and due to the necessity to verify consistency of the input quartet set.
Dense quartet sets  (such sets contain at least one quartet on every $4$-tuple of taxa) definitely suffice for this more general method
to succeed. Presented with a dense quartet set our $O(n^6)$ method will either return proof of inconsistency or determine a sparsest level-1 network
consistent with the quartet set. This answers affirmatively the question of Gambette, Berry and Paul \cite{Gambette2012} whether a level-1 phylogenetic network can be reconstructed from a dense quartet set in polynomial time.

The paper is organized as follows. In the next section, we give preliminaries on notation, cyclic orders, phylogenetic trees, and linear algebra. In Section 3 we investigate the structure of the aforementioned affine subspace of vectors which encode cyclic orders, and the relation of this subspace with phylogenetic networks. In Section 4, we describe both the restricted algorithm and the the more flexible method, and argue that the time complexity of these methods is $O(n^3)$ and $O(n^6)$, respectively.

\section{\label{sect:prelim}Preliminaries}
\subsection{Some set notation}  If $X$ is a finite set, we write
$$\binom{X}{k}:=\{Y\subseteq X\mid |Y|=k\}.$$
We abbreviate $X-x:=X\setminus \{x\}$ and $X+x:=X\cup \{x\}$.

\subsection{Cyclic orders}
If 
$S^1:=\{x\in\R^2\mid \|x\|=1\}$ is the unit circle, then an injection $\gamma: X\rightarrow S^1$ determines a {\em cyclic ordering} of $X$. The combinatorial properties of such a cyclic ordering are completely captured by the order $(x_1, x_2,\ldots, x_n)$ in which the elements of $X$ are encountered when one traverses $S^1$ in clockwise order, starting from any $s\in S^1\setminus \gamma[X]$ and continuing until all elements of $X$ are encountered. Of course, a change in the choice of $s$ will rotate the sequence.
We write $[x_1, \ldots, x_n]$ for the equivalence class of a sequence modulo the rotation, so that $[x_1, \ldots, x_n]=[x_2, \ldots, x_n, x_1]$.

If  $C=[x_1, \ldots, x_n]$ and $[y_1,\ldots, y_k]$ arises from $C$ by omitting some entries $x_i$, we say that $C$ {\em induces} $[y_1,\ldots, y_k]$ on $\{y_1,\ldots, y_k\}$, and we write $C  \succeq  [y_1,\ldots, y_k]$. We denote $C_Y :=  [y_1,\ldots, y_k]$, where $Y=\{y_1,\ldots, y_k\}$. We say that $Y\subseteq X$ is {\em consecutive} in $C$ if $Y=\{y_1,\ldots, y_k\}$ and $C=[y_1,\dots, y_k, x_1,\ldots, x_l]$, and that $C'=[y_k,\dots, y_1, x_1,\ldots, x_l]$ arises from $C$ by {\em reversing} $Y$.

If $C$ is a cyclic ordering of $X$, then the set of cyclic orderings induced on triples from $X$ characterizes $C$, see e.g. \cite{Heyting1980}.
\begin{theorem}\label{thm:circ} Let $X$ be a finite set, and let $T\subseteq\{[a,b,c]\mid a,b,c\in X\}$. Then there is a cyclic ordering $C$ such that $T=\{[a,b,c]\mid C\succeq [a,b,c]\}$ if and only if
\begin{enumerate}
\item $[a,b,c]\in T$ if and only if $[c,b,a]\not\in T$; and
\item if $[a,b,c]\in T$ and $[a,c,d]\in T$, then $[a,b,d]\in T$
\end{enumerate}
\end{theorem}
\ignore{We note that both conditions in this characterization  involve no more that four elements $a,b,c,d\in X$.
In particular, it follows that $T$ is the set of cyclic ordering if and only if each restriction of $T$ to four elements is.}

\subsection{Trees and level-1 networks}
\begin{figure}[tb]
	\subfigure[$H_1$]{
		\begin{tikzpicture}[scale=0.50,node distance=.6cm]
			\node[dot] (A) {};
			\node[scale=1.2, above left of =A] { $a$};
			\node[dot, below right of=A] (X) {};
			\node[dot, right of=X] (Y) {};
			\node[dot, right of=Y] (Z) {};
			\node[dot, above right of =Z](E){};
			\node[scale=1.2, above right of =E]{$e$};
			\node[dot, below left of=X](B){};
			\node[scale=1.2, below left of =B]{$b$};
			\node[dot, below of=Y](C){};
			\node[scale=1.2, below of =C]{$c$};
			\node[dot, below right of=Z](D){};
			\node[scale=1.2, below right of=D]{$d$};
			\draw[line] (A) -- (X) -- (Y) -- (Z) -- (E);
			\draw[line] (X) -- (B);
			\draw[line] (Y) -- (C);
			\draw[line] (Z) -- (D);
		\end{tikzpicture}
	} \hspace{.07\textwidth}
	\subfigure[$H_2$]{
		\begin{tikzpicture}[scale=.50, node distance=.6cm]
			\node[dot] (C) {};
			\node[scale=1.2, left of =C]{$c$};
			\node[dot, right of=C] (X) {};
			\node[dot, below right of=X] (Y) {};
			\node[dot, above right of=X] (Z) {};
			\node[dot, below right of=Z] (T) {};
			\node[dot, below of =Y] (D) {};
			\node[scale=1.2, below of =D]{$d$};
			\node[dot, above of =Z] (A) {};
			\node[scale=1.2, above of =A]{$a$};
			\node[dot, right of =T] (U) {};
			\node[dot, above right of =U] (B) {};
			\node[scale=1.2, above right of =B]{$b$};
			\node[dot, below right of =U] (E) {};
			\node[scale=1.2, below right of =E]{$e$};
			\draw[line] (C)--(X)--(Y)--(D);
			\draw[line] (Y)--(T)--(U)--(B);
			\draw[line] (X)--(Z)--(A);
			\draw[line] (Z)--(T);
			\draw[line] (U)--(E);			
		\end{tikzpicture}
	}\hspace{.07\textwidth}
	\subfigure[$H_3$]{
		\begin{tikzpicture}[scale=0.50, node distance=.6cm]
			
			\node[dot] (X) {};
			\node[dot, below of=X] (Y) {};
			\node[dot, below right of=Y] (Z) {};
			\node[dot, above right of=Z] (U) {};
			\node[dot, above of=U] (V) {};
			
			\node[dot, left of=X] (A) {};
			\node[scale=1.2, left of =A] { $a$};
			
			\node[dot, below left of=Y] (B) {};
			\node[scale=1.2, below left of =B] { $b$};
			
			\node[dot, below of=Z] (C) {};
			\node[scale=1.2, below of =C] { $c$};

			\node[dot, below right of=U] (D) {};
			\node[scale=1.2, below right of =D] { $d$};

			\node[dot, right of=V] (E) {};
			\node[scale=1.2, right of =E] { $e$};

			\draw[line] (X) -- (Y) -- (Z) -- (U) -- (V) -- (X);
			\draw[line] (X) -- (A);
			\draw[line] (Y) -- (B);
			\draw[line] (Z) -- (C);
			\draw[line] (U) -- (D);
			\draw[line] (V) -- (E);
		\end{tikzpicture}
	}
	\caption{\label{fig:dyadic}The three level-1 networks on 5 taxa. }
\end{figure}
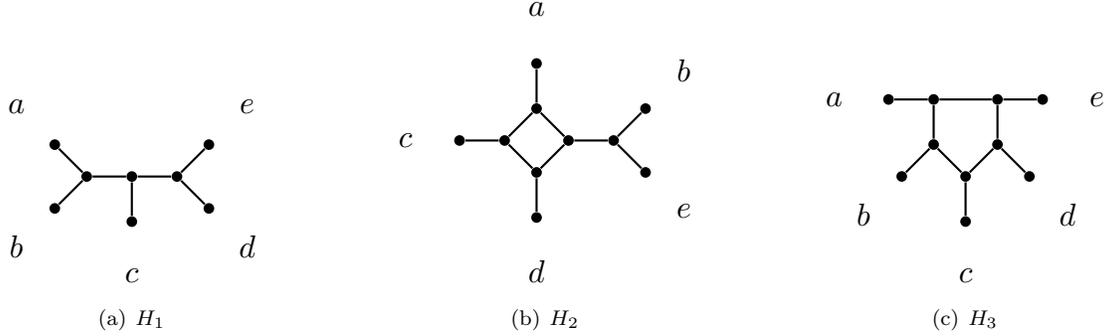

An unrooted {\em phylogenetic network on $X$\/} is an undirected simple graph $G=(V,E)$, such that $X\subseteq V$ is the set of leaves of $G$.
Such a network is {\em binary} if each vertex $v\in V\setminus X$ has degree 3. An unrooted phylogenetic network $G$ is {\em level-1\/} if   no two vertices of $G$ are connected by three internally vertex-disjoint paths.
We will simply use the term {\em level-1 network\/} in this paper to
denote an unrooted binary phylogenetic level-1 network. 
For example,there are three possible shapes for a level-1 network on five leaves, see Figure \ref{fig:dyadic}. Taking the leaf labels into account, there are several inequivalent level-1 networks for each shape. 

Given any subset $Y\subseteq X$, there is a unique subgraph  of $G$ with edges
$$F:=\{e\in E(G)\mid e \text{ is on some path between }x,y\in Y\}.$$
Deleting isolated vertices,  replacing triangles by vertices of degree 3, and suppressing vertices of degree 2, we obtain the graph $G_Y$ that $G$ {\em induces} on $Y$.

We say that a level-1 network $G$ on $X$ is {\em consistent} with a cyclic ordering $C$ if and only if there is a topological embedding of $G$ in a disk such that the restriction of the embedding to the boundary of the disk contains the elements of $X$ in the cyclic order specified by $C$. We define
$$\CC(G):=\{C\mid C\text{ a cyclic ordering consistent with }G\}.$$
We will argue below that a binary level-1 network $G$ is uniquely determined by this set $\CC(G)$. It is straightforward that if $C\in \CC(G)$ and $Y\subseteq X$, then $C_Y\in \CC(G_Y)$.

A subset $S\subseteq X$ is a {\em split} of a level-1 network on $X$ if and only if there is a cut-edge $e$ of $G$ so that $S$ equals the set of vertices from $X$ contained in one of the components of $G-e$. We say that a split is {\em nontrivial}  if $|S|, |X\setminus S|\geq 2$, and we write
$$\Sigma(G):=\{S\subseteq X\mid S\text{ is a nontrivial split of }G\}.$$
It is easy to see that if $Y\subseteq X$ and $S$ is a split of $G$ so that $|Y\cap S|, |Y\setminus S|\geq 2$, then $Y\cap S$ is a nontrivial split of $G_Y$.
We will argue below that a binary level-1 network $G$ is uniquely determined by the combination of some $C\in \CC(G)$ and the full set of nontrivial splits $\Sigma(G)$.

If $a,b,c,d$ are distinct, then the pair of pairs $\{\{ab\},\{cd\}\}$ is a  {\em quartet} on $\{a,b,c,d\}$, denoted $(ab|cd)$. If $q=(ab|cd)$ is a quartet, then the set {\em underlying} $q$ is $\underline{q}:=\{a,b,c,d\}$.
If $G$ is a level-1 network on $X$, then we say that a quartet $(ab|cd)$ is {\em compatible with $G$} if $G$ contains two disjoint paths, one from $a$ to $b$ and the other from $c$ to $d$. The set of quartets compatible with $G$ is denoted
$$\QQ(G):=\{(ab|cd)\mid (ab|cd)\text{ compatible  with } G\}.$$
Clearly, $\QQ(G_Y)=\{q\in \QQ(G)\mid \underline{q}\subseteq Y\}$ for any $Y\subseteq X$.

\subsection{\label{ss:linalg} Linear algebra} If $V$ is a vector space, and $W\subseteq V$ is an affine subspace, then we say that a vector $v\in V$ is {\em parallel to} $W$ if and only if there are $w,w'\in W$ such that $v=w'-w$. So if $W=\{x\mid Ax=b\}$, then $v$ is parallel to $W$ if $W\neq \emptyset$ and $Av=0$.

If $T$ is a finite set and $S\subseteq T$, then the {\em characteristic vector} of $S$ is $\chi^S\in \GF(2)^T$, where $\chi^S_t=1$ if and only if $t\in S$, for each $t\in T$. Finally, if $u\in \F^A$ and $B\subseteq A$, then $u_{|B}$ is the vector in  $ \F^B$ that arises by restricting $u$ to the entries in $B$.

As a part of our proposed algorithm, we solve a system of $k$ linear equations over $m$ variables over $\GF(2)$. In our application, each equation will involve at most $4$ variables and there will typically be far more equations than variables. Standard Gaussian elimination will take $O(k^2m)$ time to solve such a system.

We describe an algorithm which takes $O(km+m^3)$ time, exploiting the fact that each equation is {\em sparse} i.e. is of the form  $ax=b$ where the row vector $a$ has at most $O(1)$ nonzero entries. We below assume that such $a$ are stored by listing their nonzero entries, so that computing the inner product $ax$ takes $O(1)$ time for any $x$.

Consider the following auxiliary problem:
\begin{tabbing}
{\bf Given:} \= A sparse linear equation $ax=b$, a vector $u$, and a $m\times d$ matrix $V$ with independent columns.\\
{\bf Find:}\> Either the signal that $U'=\emptyset$ or a vector $u'$ and a matrix $V'$ with $d'$ independent columns such \\
that $U'=\{u'+V'y'\mid y'\in \GF(2)^{d'}\}$,
where $U':=\{x\in U\mid ax=b\}$ and  $U:=\{u+Vy\mid y\in \GF(2)^d\}$.
\end{tabbing}
The standard procedure for solving this problem takes $O(d)$ to compute $aV$, using the sparsity of $a$. If $aV=0$, then $U'=U$ or $U'=\emptyset$ depending on whether $au=b$ so that we are the done in $O(d)$ overall time. Otherwise, it takes $O(dm)$ time to perform the necessary column operations to compute $u', V'$ from $u,V$.

We return to the main problem:
\begin{tabbing}
{\bf Given:} \= Sparse linear equations $a_1x=b_1$,\ldots, $a_kx=b_k$ over $\GF(2)$ in $m$ variables.\\
{\bf Find:}\> Either the signal that $U=\emptyset$, or a vector $u$ and a matrix $V$ with $d$ independent columns such\\
that $U=\{u+Vy\mid  y\in \GF(2)^{d}\}$, a set $K$ such that $U=\{x\in \GF(2)\mid a_ix=b_i, i\in K\}$ with $|K|+d=m$, \\
where $U:=\{x\in \GF(2)^m\mid a_1x=b_1,\ldots, a_kx=b_k\}$.
\end{tabbing}

For $j=0,1,2,\ldots$ we iteratively compute $K_j$, a vector $u^j$, and a matrix $V^j$ with $d_j$ independent columns, so that
$$U_j:=\{x\mid a_ix=b_i\text{ for }i=1,\ldots, j\}=\{x\mid a_ix=b_i\text{ for }i\in K_j\}=\{u^j+V^jy\mid y\in \GF(2)^{d_j}\}.$$
The algorithm starts with $U_0=\GF(2)^m$, $u^0=0, V=I$, $K_0=\emptyset$, $j=0$.
In the $j$-step, we compute $u^{j+1}, V^{j+1}, d_{j+1}$ from $u^j, V^j, d_j$, which amounts to solving the above auxiliary problem.
If it turns out that  $U_{j+1}=\emptyset$, then we terminate with the signal that $U=\emptyset$.
If  $d_{j+1}=d_j$, we put $K_{j+1}=K_j$, and if $d_{j+1}=d_j-1$, then we set $K_{j+1}=K_j\cup\{j+1\}$.
Then we put $j\leftarrow j+1$. We proceed until $j=k$ if possible, and put $u=u^j$, $V=V^j$, $d=d_j$, $K=K_j$ at termination.

In this algorithm, the $j$-th stage takes $O(m^2)$ if $j+1\in K$, and $O(m)$ time otherwise. As $|K|+d=m$, we have $O(m)$ stages taking $O(m^2)$ time and $O(k)$ stages taking $O(m)$ time. All together $O(km+m^3)$ time.

\section{\label{sect:space}Level-1 networks, cyclic orderings and an affine space}
\subsection{Encoding cyclic orderings} With any cyclic ordering $C$ of $X$, we associate a vector $u^C\in \GF(2)^{X^3}$ by setting
$$u^C(x,y,z):=\left\{\begin{array}{ll} 1  &\text{if }x,y,z\text{ are distinct and }C\succeq [x,y,z]\\ 0&\text{otherwise}\end{array}\right.$$
for each triple of taxa $x,y,z\in X$. We say that a vector $u\in \GF(2)^{X^3}$ is {\em cyclic} if $u=u^C$ for some cyclic ordering $C$ of $X$. We say that a subspace of $ \GF(2)^{X^3}$ is cyclic if all the vectors it contains are cyclic.

For ease of notation, we define the following. If $u\in \GF(2)^{X^3}$ and $Y\subseteq X$, then $u_Y:=u_{|Y^3}\in \GF(2)^{Y^3}$  denotes the restriction of $u$ to the entries in $Y^3$.  For a subset $U\subseteq \GF(2)^{X^3}$, we put $U_Y:=\{u_Y\mid u\in U\}$.
\begin{lemma}\label{lemma:cyclic} Let $u\in \GF(2)^{X^3}$ be a cyclic vector. Then
\begin{enumerate}
\item[(0)] $u(x,x,y)=0$ for all $x,y\in X$
\item[(1)] $u(x,y,z)+u(y,z,x)=0$ for all $x,y,z\in X$
\item[(2)] $u(x,y,z)+u(y,x,z)=1$ for all $x,y,z\in X$
\item[(3)] $u(t,x,y)+u(t,x,z)+u(t,y,z)+u(x,y,z)=0$ for all $t,x,y,z\in X$
\end{enumerate}
\end{lemma}
\proof $(1)$ and $(2)$ follow from the definition of $u^C$ and from the fact that $[x,y,z]=[y,z,x]$, and $(2)$ is immediate from Theorem  \ref{thm:circ}. To see $(3)$, note that if $u$ is cyclic, then $u_Z$ is cyclic for $Z=\{t,x,y,z\}$.  Inspection of the six cyclic vectors in $\GF(2)^{Z^3}$  (see Table \ref{tab:four}) then proves the statement.\endproof
\begin{table}
\begin{tabular}{c|cccc}
 $C$& $u^C(t,x,y)$& $u^C(t,x,z)$ &$u^C(t,y,z)$ & $u^C(x,y,z)$\\ \hline
$[t,x,y,z]$ &1 & 1  &  1& 1 \\
 $[t,x,z,y]$& 1 & 1 & 0 & 0  \\
$ [t,y,x,z$]& 0 & 1 & 1 & 0   \\
$[t,y,z,x]$& 0 & 0 & 1 & 1  \\
$[t,z,x,y]$& 1 & 0 & 0 & 1  \\
$[t,z,y,x]$& 0 & 0 & 0 & 0  \\
   \end{tabular}

\caption{\label{tab:four} The six cyclic vectors in $\GF(2)^{Z^3}$, where $Z={\{t,x,y,z\}}$}\end{table}

We define the affine subspace $\UU^X\subseteq \GF(2)^{X^3}$ as the set of vectors $u$ satisfying the conditions $(0), (1), (2), (3)$ of Lemma \ref{lemma:cyclic}. It is straightforward that if $u\in \UU^X$ and $Y\subseteq X$, then $u_Y\in \UU^Y$.
\begin{lemma}\label{lemma:quad} Let $u\in\UU^X$. Then $u$ is cyclic if and only if
\begin{equation}\label{ucirc}u(t,x,y)u(t,y,z)+u(t,x,z)u(x,y,z)=0\text{ for all }t,x,y,z\in X.\end{equation}
\end{lemma}
\proof Let $T:=\{[x,y,z]\mid u(x,y,z)=1\}$. As $u\in \UU^X$, the set $T$ satisfies the first condition in Theorem  \ref{thm:circ}.
The second condition on $T$ is equivalent to \eqref{ucirc}.
So \eqref{ucirc} holds for $u$ if and only $T=\{[x,y,z]\mid C\succeq [x,y,z]\}$ for some cyclic ordering $C$ of $X$ if and only if  $u=u^C$ for some cyclic ordering $C$ of $X$.
\endproof
\begin{lemma} \label{lemma:four} Let $u\in \UU^X$. 
Then $u$ is 
cyclic if and only if 
$u_Z$ is cyclic for each $Z\in\binom{X}{4}$. 
\end{lemma}
\proof Immediate from  Lemma \ref{lemma:quad}.\endproof

In the remainder of the section, we fix the set of taxa $X$ and write $\UU=\UU^X$. By  Lemma \ref{lemma:cyclic}, $\UU$ contains all cyclic vectors $u\in \GF(2)^{X^3}$. Note that the six cyclic vectors of  Table \ref{tab:four} are exactly the vectors that satisfy \eqref{ucirc}. The two other vectors of $\UU^Z$, a space of dimension 3 over $\GF(2)$ with 8 elements, are not cyclic.

\subsection{An affine space from a level-1 network} For a level-1 network $G$ on $X$, we define
$$U(G):=\{u^C\mid C\in \CC(G)\}.$$
We can explicitly describe $U(G)$ in terms of the splits of $G$.
For a subset $Y\subseteq X$, we define a vector $v^Y\subseteq \GF(2)^{X^3}$ by setting
$$v^Y(x,y,z):=\left\{\begin{array}{ll} 1  &\text{if }x,y,z\text{ are distinct and }|Y\cap\{x,y,z\}|\geq 2\\ 0&\text{otherwise}\end{array}\right.$$

\begin{lemma} \label{lemma:char}Let $G$ be a binary level-1 network on $X$, and let  $Y\subseteq X$. Then the following are equivalent:
\begin{enumerate}
\item $Y$ is a split of $G$;
\item $u+v^Y\in U(G)$ for all $u\in U(G)$; and
\item  $u+v^Y\in U(G)$ for some $u\in U(G)$.
\end{enumerate}
\end{lemma}
\proof
$(1)\Rightarrow (2)$: Suppose $Y$ is a split of $G$.  Let $u\in U(G)$ be any vector, and let $C\in\CC(G)$ be such that  $u=u^C$. Since $C$ is compatible with $G$, $Y$ is consecutive in $C$. Suppose $C'$ is the cyclic ordering that arises from $C$ by reversing the order of $Y$ in $C$. Then $C'$ is again compatible with $G$, and  $u+v^Y=u^C+v^Y=u^{C'}\in U(G)$.

$(2)\Rightarrow (3)$: Trivial.

$(3)\Rightarrow (1)$: Let $C$ be such that $u'=u^C+v^Y\in U(G)$. Then $Y$ is consecutive in $C$. If not, there are $a,c\in Y$ and $b,d\not\in Y$ so that $[a,b,c,d]\preceq C$. Then $u'(a,b,c)=u'(a,c,d)=0$ and $u'(a,b,d)=u'(b,c,d)=1$, contradicting that $u'$ is cyclic.
So we have $u'=u^{C'}$, where $C'$ arises by reversing $Y$ in $C$. Since $C, C'\in \CC(G)$, there exist outerplanar drawings of $G$ that exhibit $C$ and $C'$ on the boundary, and so there are no two vertex-disjoint paths from $Y$ to $X\setminus Y$, for the reversal of $Y$ makes that in one of the drawings the paths would cross. Since the vertices of $V(G)\setminus X$ have degree 3, this implies that $G$ does not have two edge-disjoint paths from $Y$ to $X\setminus Y$. By Menger's Theorem, there is a set $W\subseteq V(G)$ so that $Y\subseteq W$, $(X\setminus Y)\cap W=\emptyset$, and $|\delta_G(W)|\leq 1$. Hence $Y$ is a split  of $G$, as required.
\endproof
It follows that given any fixed $u\in U(G)$ we may determine $\Sigma(G)$ as
$$\Sigma(G)=\{Y\subseteq X\mid |Y|, |X\setminus Y|\geq 2, u+v^Y\in U(G)\} $$
Moreover, any $u\in U(G)$ determines a $C\in\CC(G)$ so that $u=u^C$, so that it will be possible to reconstruct $G$ from a description of $U(G)$. We will explain how to do this efficiently in section \ref{sect:algorithm}. The next theorem  shows that in turn, $U(G)$ can be described in terms of $\Sigma(G)$ and some $C\in\CC(G)$.
\begin{theorem} \label{thm:char} Let $G$ be a binary level-1 network on $X$, and let $C$ be some cyclic ordering of $X$ compatible with $G$. Then
$$U(G)=\{u^C+\sum_{Y\in\Sigma(G)} \alpha_Yv^Y\mid \alpha_Y\in \GF(2) \} .$$
\end{theorem}
\proof $\supseteq$: Let $u=u^C+\sum_{Y\in S} v^Y$, for some $S\subseteq\Sigma(G)$. Then $u\in U(G)$ by induction on $|S|$ and using  Lemma \ref{lemma:char} $(1)\Rightarrow (2)$ in the induction step.

$\subseteq$: suppose that $u\in U(G)$, say $u=u^{C'}$. Since $C$ and $C'$ are both compatible with $G$, there is a sequence of orderings $C=C_1,\ldots, C_k=C'$ so that
$C_{i+1}$ arises from $C_i$ by reversing $Y_i=X\cap W_i$, where $W_i$ is a component of $G\setminus e_i$  for some $e_i\in E(G)$. Then  $u^{C_i}=u^{C_{i+1}}+v^{Y_i}$. The theorem follows.
\endproof
Since $U(G)$ contains only cyclic vectors, we clearly have  $U(G)\subseteq \UU$, and it is immediate from the theorem that $U(G)$ is an affine subspace of $\UU$.

\subsection{An affine space from quartets} If we are given a set of quartets on $X$, an affine subspace of $\UU$ can be constructed from it, since every quartet gives rise to one equation over $GF(2)$.
\begin{lemma} \label{lemma:char2} Let $G$ be a binary level-1 network on $X$, and let $a,b,c,d\in X$.
Then $(ab|cd)$ is a quartet compatible with $G$ if and only if
$$u(a,b,c)+u(a,b,d)=0$$
for all $u\in U(G)$.
\end{lemma}
\proof If $(ab|cd)\in\QQ(G)$, then for any $C\in \CC(G)$, we have one of $C\succeq [a,b,c,d], [a,b,d,c], [b,a,c,d]$, or $[b,a,d,c]$. In any case, we have $u^C(a,b,c)+u^C(a,b,d)=0$.

Conversely, suppose that
$(ab|cd)\notin \QQ(G)$. Then either $(ac|bd)\in \QQ(G)$ or $(ad|bc)\in \QQ(G)$ or both. 
If $(ac|bd)\in \QQ(G)$, but $(ad|bc)\notin \QQ(G)$ then there is a split $S\in \Sigma(G)$ with $a,c\in S$ but $b,d\notin S$. Theorem \ref{thm:char} implies that a $u\in U(G)$ exists such that $u(abc)\neq u(abd)$. Similarly
if $(ad|bc)\in \QQ(G)$ but $(ac|bd)\notin \QQ(G)$. If both quartets are in $\QQ(G)$, then $u(a,c,b)=u(a,c,d)$ and $u(a,d,b)=u(a,d,c)$  for all
$u\in U(G)$ implying that $u(a,b,c)=u(a,d,c)=u(a,d,b)\neq u(a,b,d)$ for all $u\in U(G)$.
\endproof

For a set of quartets $Q$, we write
$$U(Q):=\{u\in\UU\mid u(a,b,c)+u(a,b,d)=0\text{ for all }(ab|cd)\in Q\}.$$

We give a second characterization of  $U(G)$, in terms of the quartets of $G$.
\begin{theorem} \label{thm:char2} Let $G$ be a binary level-1 network on $X$. Then $U(G)= U(\QQ(G)).$
\end{theorem}
\proof Let $X, G$ be a counterexample with $|X|$ as small as possible. The case $|X|= 4$ is straightforward to verify as there are only two possible binary level-1 networks on any set of taxa  $X$ of that size, up to isomorphism. By Lemma \ref{lemma:char2}, we have $U(G)\subseteq U(\QQ(G))$. As $G$ is a assumed to be a counterexample, there must exists $u\in U(\QQ(G))\setminus U(G)$. There are two cases to consider: that $u$ is not cyclic, and that $u$ is cyclic.

If $u$ is not cyclic, then there is some 4-tuple $Z$ from $X$ so that $u_Z$ is not cyclic, by Lemma \ref{lemma:four}. Then $u_Z\not\in U(G_Z)$ but $u_Z\in U(\QQ(G))_Z\subseteq U(\QQ(G_Z))$. Then by minimality of $|X|$, we have $X=Z$, i.e. $|X|=4$.

If on the other hand $u$ is cyclic, then $u=u^C$ for some cyclic ordering $C$ of $X$. Pick an $x\in X$,  let $y,z\in X$ be the two vertices adjacent to $x$ in $C$, and let $C':=C_{X- x}$.  As $X$ was a minimal counterexample, we have $u^{C'}=u_{X- x}\in U(\QQ(G_{X- x}))=U(G_{X- x})$, so that $C'$ is compatible with $G_{X- x}$.
Consider an outerplanar drawing $f$ of $G_{X- x}$  that realizes $C'$ on the boundary of a disk $D$. Then the face of the drawing incident with $y$ and $z$ is bounded by the segment $yz$ on the boundary of $D$ and a path $P$ from $y$ to $z$ in $G_{X-x}$.  Since $u^C\in U(\QQ(G))$, any path $Q$ from $x$ to a vertex from $X\setminus\{y,z\}$ in $G$ must intersect $P$.

The graph $G_{X- x}$ arises from $G$ by both deleting vertex $x$ and the unique  edge $xs$ incident with $x$, and either
\begin{enumerate}
\item eliminating a possible circuit $(s,t,u,w)$ through $s$ by deleting $s$ and suppressing the vertices  $t,w$, or
\item suppressing the degree-2 vertex $s$, if $s$ is not on such a cycle.
\end{enumerate}
In the latter case, it is possible to extend $f$ to a drawing of $G$ which realizes $C$ on the boundary of the disk $D$ if and only if $s$ lies on the path $P$, by drawing the edge from $x$ to $s$ inside $F$. For if $s$ does not lie on $P$, it is possible to construct a path $Q$ disjoint from $P$ to some vertex in $X\setminus\{y,z\}$, which would give a contradiction. In the former case, we can similarly  extend $f$ if and only if both $t$ and $w$ lie on $P$. Again, it is possible to construct a path $Q$ from $x$ to $X\setminus\{y,z\}$ if one of $t,w$ is not on $P$, contradiction.\endproof
We will argue in the next section that there is a considerable degree of redundancy in the system of $O(n^4)$ equations that determine $U(\QQ(G))$ in the above theorem, and we will show that a certain set of $O(n^3)$ equations will also do.

A set of quartets $Q$ is {\em dense\/} if it contains a quartet on every four taxa, i.e. 
if for any $Z\in \binom{X}{4}$ there is a $q\in Q$ such that $\underline{q}=Z$.
The following theorem is a direct consequence of Lemma \ref{lemma:cyclic}, Lemma \ref{lemma:quad} and Lemma \ref{lemma:char2}.
\begin{theorem}\label{theorem:dense}
If $Q$ is a dense quartet set, then $U(Q)$ is cyclic, and can be characterized as
\begin{eqnarray*}U(Q)=\{ u\in GF(2)^{X^3}&\mid& u(a,b,c)=u(b,c,a)=u(a,c,b)+1  \mbox{ for all } a,b,c \in X, \\
&&  u(a,b,c)=u(a,b,d) \mbox{ and } u(a,c,d)=u(b,c,d) \mbox{ for all } (ab|cd)\in Q\}.
\end{eqnarray*}
 \endproof
\end{theorem}

Splits of $G$ give rise to decompositions of $U(G)$, which will be useful for proving, at the end of this subsection, that
$U(G)$ completely determines $G$ if $G$ is a level-1 network.  
First, let us define sums and decompositions of spaces.

Let $X_1$, and $X_2$ be finite sets, 
with a single element $z$ in their intersection: $X_1\cap X_2=\{z\}$. Let $X=(X_1\cup X_2)-z$. If $u_1\in \UU^{X_1}$, and $u_2\in \UU^{X_2}$, then the {\em sum\/} $u=u_1\oplus u_2$ of these two vectors is the vector $u\in \UU^X$ that satisfies
\[u(a,b,c)=
\left\{ \begin{array}{ll}
u_1(a,b,c) & \mbox{ if } a,b,c\in X_1 \\
u_1(a,b,z)& \mbox{ if } a, b\in X_1, c\in X_2 \\
u_2(a,b,z) &\mbox{ if } a, b\in X_2, c\in X_1 \\
u_2(a,b,c) &\mbox{ if } a, b, c\in X_2. \end{array} \right.
\]
 The {\em sum\/} $U=U_1\oplus  U_2$ of two subspaces $U_1\subseteq \UU^{X_1}$ and $U_2\subseteq \UU^{X_2}$ is then defined  as the following subspace of $\UU^X$
\[U_1\oplus U_2:=\{u_1\oplus u_2\mid u_1\in U_1, \; u_2\in U_2\}.\]

Inspired by Lemma \ref{lemma:char}, where splits of a level-1 network $G$ are characterized as subsets $Y$ of the set of taxa $X$ such that $v^Y$ is parallel to the subspace $U(G)$, we define the set of (nontrivial) \begin{em} splits\/ \end{em} of an arbitrary subspace $U$ of $U^{X}$ as
\[\Sigma(U):=\{S\subseteq X\mid |S|\geq 2, \; |X\setminus S|\geq 2, \; v^S \mbox{ is parallel to } U\}=\]\[\{S\subseteq X\mid |S|\geq 2, \; |X\setminus S|\geq 2, \; v^S+u\in U \mbox{ for all } u\in U\} \]

Note that a sum $U=U_1\oplus U_2$ of two subspaces $U_1\subseteq \UU^{X_1}$ and $U_2\subseteq \UU^{X_2}$ with $X_1\cap X_2=\{z\}$ as defined above always has a nontrivial split $X_1-z$  if $1_{X_1-z}$ is parallel to $U_1$, and if $|X_1-z|, |X_2-z|\geq 2$.
Conversely, cyclic quartet spaces having a nontrivial split can be decomposed (written as a sum).
\begin{theorem}\label{thm:decomp}
Let $Q$ be a set of quartets on $X$, such that $U(Q)$ is cyclic. If $S\in \Sigma(U(Q))$,
then there are sets of quartets $Q_1$ on $S\cup\{z\}$ and $Q_2$ on $(X\setminus S)\cup\{z\}$, where $z$ is a new taxon, such that
\[U(Q)=U(Q_1)\oplus U(Q_2).\] Moreover,
if $Q=\QQ(G)$ and $S$ is a nontrivial split of $G$ (associated with an edge $e$) for some connected level-1 network $G$,  then $Q_i=\QQ(G_i)$, for $i=1,2$, where $G_1$ is the level-1 network obtained from $G$
by contracting one component of $G-e$ to 
a single leaf $z$,  and $G_2$ is the  level-1 network obtained from $G$ by contracting
the other component of $G-e$ to a single leaf $z$.
\end{theorem}
\proof
Without loss of generality, $Q$ is maximal in the sense that for all quartet sets $Q'$ such that $U(Q')=U(Q)$, it holds that $Q'\subseteq Q$. Then in particular, $Q$ is dense.
Since $v^S+u\in U(Q)$ for any $u \in U(Q)$,
\[v^S(a,b,c)=v^S(a,b,d) \mbox{ and } v^S(a,c,d)=v^S(b,c,d) \mbox{ for all } (ab|cd)\in Q.\]
This implies that for any $a,b\in S$ and $c,d\in X\setminus S$ it holds that
$u(a,b,c)=u(a,b,d)$ for all $u\in U(Q)$. Indeed, if there would be an $u\in U(Q)$ with $u(a,b,c)\neq u(a,b,d)$, then it can be checked in Table \ref{tab:four} that either $u$ is not cyclic, or $v^S+u$ is not cyclic, contradicting the fact that $U(Q)$ is cyclic and $v^S$ is parallel to $U(Q)$.
By maximality of $Q$, $(ab|cd)\in Q$ for any $a,b\in S$ and $c,d\in X\setminus S$. Moreover, $(ab|cd)$ is the only
quartet in $Q$ on such a quadruple $a,b,c,d$ with $a,b\in S$ and $c,d\in X\setminus S$ (because $v^S(a,c,b)\neq v^S(a,c,d)$ and $v^S(a,d,b)\neq v^s(a,d,c)$) by definition of $v^S$).

Now choose $c_1\in X\setminus S$ and $c_2\in S$ arbitrarily, and define (dense) quartet sets $Q_1$ on $S\cup\{z\}$ and $Q_2$ on $(X\setminus S)\cup\{z\}$, where  $z$ is a new taxon, by
\[Q_1:=\{(ab|zd)\mid (ab|c_1d)\in Q, a,b,d\in S\}\cup\{(ab|cd)\mid (ab|cd)\in Q,  a,b,c,d\in S\}\]
\[Q_2:=\{(ab|zd)\mid (ab|c_2d)\in Q, a,b,d\notin S\}\cup\{(ab|cd)\mid (ab|cd)\in Q,  a,b,c,d\notin S\}\]
Then
$U(Q)=U(Q_1)\oplus U(Q_2)$.
Indeed, if  $u\in U(Q)$, then consider the restrictions $u_1'=u_{S+c_1}$ and $u_2'=u_{(X\setminus S)+c_2}$ of $u$ to the subsets $S+c_1$ and $X\setminus S+c_2$ respectively. Define $u_1\in \UU^{S+z}$ by $u_1(a,b,c)=u_1'(a,b,c)$ if $a,b,c\in S$ and $u_1(a,b,z)=u(a,b,c_1)$ if $a,b\in S$. Define $u_2$ similarly. Then it is easy to check that $u=u_1\oplus u_2$
(using that $u(a,b,d)=u(a,b,c_1)$ for all $a,b\in S$ and $d\notin S$).
Conversely, let $u_1\in U(Q_1)$ and $u_2\in U(Q_2)$, and let $u=u_1\oplus u_2$. We show that $u\in U(Q)$ by proving that $u(a,b,c)=u(a,b,d)$ for all $(ab|cd)\in Q$ (using the characterizations of $U(Q_1), U(Q_2)$, and $U(Q)$ from Theorem \ref{theorem:dense}). So let $(ab|cd)\in Q$. If $a,b\in S$, $c,d\notin S$, then $u(a,b,c)=u_1(a,b,z)=u(a,b,d)$ by definition. If $a,b,d\in S$, $c\notin S$, then because $(ab|cc_1), (ab|cd) \in Q$, we have for all $x\in U(Q)$ that $x(a,b,c_1)=x(a,b,c)=x(a,b,d)$ and by maximality of $Q$ that $(ab|c_1d)\in Q$ so $(ab|zd)\in Q_1$ and hence $u(a,b,c)=u_1(a,b,z)=u_1(a,b,d)=u(a,b,d)$. The other possibilities are symmetrical. This proves that $u\in U(Q)$.

Finally, suppose that $Q=\QQ(G)$ for a connected level-1 network $G$, and $S$ is a (nontrivial) split of $G$. Then there is an edge $e$ in $G$, such that $G-e$ has two components $H_1$ and $H_2$, and $S$ is the set of leaves contained in the component $H_1$. Now it is not hard to show that $Q_1=\QQ(G_1)$ for the graph $G_1$ obtained from $G$ by contracting the vertices of $H_2$ to a single leaf $z$ (attached to the edge $e$). 
Similarly $Q_2=\QQ(G_2)$ for the graph $G_2$ obtained from $G$ by contracting the vertices of $H_1$ to a single leaf $z$.  
\endproof
\begin{theorem}\label{thm:unique}
Let $G$ and $G'$ be connected level-1 networks on $X$ such that $U(G)=U(G')$. Then $G=G'$.
\end{theorem}
\proof
Let $G$, $G'$ form a counterexample with $|X|$ minimal.
By Lemma \ref{lemma:char}, $S$ is a nontrivial split of $G$ if and only if it is a nontrivial split of $U(G)=U(G')$
if and only if it is a nontrivial split of $G'$.
By Theorem \ref{thm:char2}, $U(G)=U(\QQ(G))$ and $U(G')=U(\QQ(G'))$.
If $G$ has a nontrivial split $S$, then $U(G)=U(G')$ decomposes by Theorem \ref{thm:decomp} as
\[U(G)=U(G')=U(Q_1)\oplus U(Q_2), \] where $Q_1$ is the quartet set of a network $G_1$ arising from  $G=(V,E)$ by
contracting the vertex set  $W\subseteq V$ of one component of $G-e$ (with $W\cap X=S$) for some edge $e\in E$ to a single vertex in $G$, and it is also the quartet set of a network $G_1'$ arising from $G'=(V',E')$ by contracting some vertex set $W'\subseteq V'$ of one component of $G'-e'$ (with $W'\cap X=S$) for some edge $e'\in E'$ to a single vertex in $G'$. So $U(G_1)=U(Q_1)=U(G_1')$ by Theorem \ref{thm:char2}, and by minimality of the counterexample $G_1=G_1'$. Similarly $U(G_2)=U(Q_2)=U(G_2')$ for networks $G_2$ and $G_2'$ arising from $G$ and $G'$ respectively by contracting the vertex sets $V\setminus W$ in $G$ and $V'\setminus W'$ in $G'$, respectively (each with leaf set $X\setminus S$), and again by  minimality of the counterexample we obtain $G_2=G_2'$.
But then $G=G'$ because $G$ can be uniquely recovered from $G_1$ and $G_2$.

If $G$ has no nontrivial splits, then it is a single vertex, an edge, a star, or a circuit with leaves, and the same is true for $G'$.
 The theorem follows from the fact that two such level-1 networks are isomorphic if they have the same leaf set and allow the same cyclic orderings of their leaves.
\endproof

\subsection{Witnesses}
By Theorem \ref{thm:char2}, we know that if $u\in \UU\setminus U(G)$, then $u\in \UU\setminus U(\QQ(G))$, so that for some $(ab|cd)\in \QQ(G)$, we then have  $u(a,b,c)+u(a,b,d)\neq 0$. This is a local condition on $u$, in the sense that the equation only involves entries indexed by triples from $\{a,b,c,d\}^3$. So the quadruple $\{a,b,c,d\}$ can be considered a witness for the fact that $u\not\in U(G)$. We will argue that such witnesses are not isolated.

We investigate the local witnesses for $u\not \in U(G)$, where $G$ is a level-1 network. First, we consider the situation that $u$ is not even cyclic.
\begin{lemma} \label{lemma:noncyclic}Let $Y$ be a set with 5 elements, and let $u\in \UU^Y$. If $u$ is not cyclic, then there are exactly two $Z\in\binom{Y}{4}$ so that $u_Z$ is not cyclic.\end{lemma}
\proof Let $Y=\{a,b,c,d,e\}$. We first argue that the number of $Z\in\binom{Y}{4}$ so that $u_Z$ is not cyclic is even.
For each $Z=\{t,x,y,z\}$, Lemma \ref{lemma:quad} states that $\alpha(t,x,y,z)=1$ if and only if $u_Z$ is not cyclic, where $\alpha(t,x,y,z):=u(t,x,y)u(t,y,z)+u(t,x,z)u(x,y,z)$.
Consider the sum $\alpha:=\alpha(a,b,c,d)+\alpha(a,b,c,e)+\alpha(a,b,d,e)+\alpha(a,c,d,e)+\alpha(b,c,d,e)$.
It is straightforward that $\alpha=0$, using the equations in Lemma \ref{lemma:cyclic} to rewrite $\alpha$ to an expression involving only $u(x,y,z)$ with $x=a$ and $y<z$ where say, $b<c<d<e$. Hence, $u_Z$ is not cyclic for an even number of $Z\in\binom{Y}{4}$.
It remains to show that there is no $u\in \UU^Y$ such that $u_Z$ is not cyclic for exactly four $Z\in\binom{Y}{4}$.
By symmetry, we may assume that $u_Z=u^{[a,b,c,d]}$ if $Z=\{a,b,c,d\}$ and that $u_Z$ is not cyclic for the other 4-sets $Z$ from $Y$.
A straightforward case-check shows that $u^{[a,b,c,d]}$ cannot be extended to such a vector in $\UU^Y$.\endproof

Next, we argue a similar condition for the 4-sets of taxa which witness that an $u\in \UU$ is not in $U(G)$.
\begin{lemma} \label{lemma:five}Let $Y$ be a set with 5 elements, let $G$ be a level-1 network on $Y$, and let $u\in \UU^Y\setminus U(G)$. Then there are at least two $Z\in\binom{Y}{4}$ so that $u_Z\not \in U(G)_Z$.\end{lemma}
\proof Let $u\in \UU^Y\setminus U(G)$. If $u$ is not cyclic, then we are done by Lemma \ref{lemma:noncyclic}. So $u=u^C$ for some cyclic ordering of $Y$. We may assume that $Y=\{a,b,c,d,e\}$ and using symmetry, that either $G=H_1, H_2$, or $H_3$, where $H_1, H_2, H_3$ are as in Figure \ref{fig:dyadic}. For each $G$, there are at most 24 eligible cyclic orderings $C$ of $Y$. We omit the details of the remaining finite case-check.\endproof

\begin{theorem}\label{thm:rooted} Let $G$ be a level-1 network on $X$, let $\infty\in X$ and let $Q=\{q\in \QQ(G)\mid \infty\in \underline{q}\}$. Then $U(Q)=U(G)$.\end{theorem}
\proof Suppose not. Then there is some $u\in U(Q)\setminus U(G)$. As $U(G)=U(\QQ(G))$ by Theorem \ref{thm:char2}, there is some  $(ab|cd)\in \QQ(G)$ so that $u(a,b,c)+u(a,b,d)\neq 0$. Then $\infty\not\in \{a,b,c,d\}$ by definition of $Q$ and $U(Q)$. Applying Lemma \ref{lemma:five} to the 5-tuple $Y=\{a,b,c,d,\infty\}$, we arrive at a contradiction, since  $u_Z\in U(G)_Z$ for exactly four of the sets $Z\in\binom{Y}{4}$.
\endproof

\subsection{\label{ss:basis} A basis for $\UU$}
In what follows, we assume that the set of taxa $X$ contains a fixed taxon $\infty$, and we fix a linear ordering $<$ of $X$.
We will consider a parametrization of the affine space $\UU$ relative to $\infty$ and $<$.
Consider the set of triples $T:=\{(\infty, a,b) \mid a<b, ~a,b\in X\}$. The following is straightforward from Lemma \ref{lemma:cyclic}.
\begin{lemma} For any $y\in \GF(2)^T$, there is a unique vector $u\in \UU$ so that $u_{|T}=y$.\end{lemma}
So for any affine subspace $U\subseteq \UU$ we have $\#U=\#U_{|T}$, so that $\dim(U)=\dim(U_{|T})$. In particular, we have $\dim(\UU)=|T|=\binom{n-1}{2}$, where $n=|X|$.

If $a=\infty$  then the linear equation in Lemma \ref{lemma:char2} refers only entries $u(\infty, x,y)$ of $u$. Using that $u(\infty,x,y)=u(\infty,y,x)+1$ if necessary, the equation can be seen to be equivalent to one on $u_{|T}$.
If $\infty\not\in\{a,b,c,d\}$ then
$$u(a,b,c)+u(a,b,d)=u(\infty,a,c)+u(\infty,b,c)+u(\infty,a,d)+u(\infty,b,d).$$
so that again the linear equation of Lemma \ref{lemma:char2} may be rewitten to an equivalent equation on $u_{|T}$.

We conclude that given a set of quartets $Q$, computing $U(Q)_{|T}$ amounts to finding the set of solutions to a system of $|Q|$ linear equations in $\binom{n-1}{2}$ variables over $\GF(2)$.

\subsection{The dimension of cyclic spaces}
Our main example of an affine subspace $U\subseteq \UU$ such that each vector in $U$ is cyclic is the space $U(G)$ derived from a level-1 network $G$. We bound the dimension of such spaces in the next lemma. In what follows, ${\bf 1}_X$ denotes the vector from $\GF(2)^{X^3}$ such that ${\bf 1}_X(x,y,z)=1$ if and only $x,y,z\in X$ are distinct. Note that ${\bf 1}_X$ is parallel to $U(Q)$ for any set of quartets $Q$ on $X$, since each of the equations defining $U(Q)$ involves an even number of variables. Also, observe that if $C$ is a cyclic ordering of $X$ and $C'$ arises by reversing $C$, then $u^C+u^{C'}={\bf 1}_X$.
\begin{lemma} \label{lemma:dimG} If $G$ is a binary level-1 network on $X$, then $\dim(U(G))=|\Sigma(G)|/2 +1$. Consequently,
$\dim(U(G)) \leq |X|-2$, with equality if and only if $G$ is a tree.\end{lemma}
\proof By Theorem \ref{thm:char}, we have $U(G)=\{u+\sum_{Y\in \Sigma(G)} \lambda_Y v^Y\mid \lambda_Y\}$. The dimension of $U(G)$ therefore equals the maximum number of linearly independent vectors among $\{v^Y\mid Y\in \Sigma(G)\}$. Since $v^Y+v^{X\setminus Y}={\bf 1}_X$, each of these vectors is spanned by
$$\{{\bf 1}_X\}\cup \{v^Y\mid Y\in \Sigma(G), \infty\in Y\},$$
where $\infty\in X$ is arbitrarily chosen.  It is straightforward that the latter set of vectors is linearly independent, so the dimension of $U(G)$ equals the cardinality of this set of vectors. The elements of $\{Y\in \Sigma(G)\mid \infty \in Y\}$ correspond 1-1 to cut-edges of $G$ not incident with $X$. There are at most $|X|-3$ such edges in a binary level-1 network on $X$, with equality if and only if $G$ is a tree.
\endproof
The same conclusion about the dimension follows from a much weaker condition.
\begin{theorem} \label{thm:dim} Let $U\subseteq \UU^X$ be an affine subspace so that ${\bf 1}_X$ is parallel to $U$. If each vector $u\in U$ is cyclic, then $\dim(U)\leq |X|-2$, with equality being attained if and only if $U=U(G)$ for some binary tree $G$ on $X$.
\end{theorem}
\proof The theorem is straightforward if $|X|\leq 4$. We first show that $\dim(U)\leq |X|-2$ by induction on $|X|$.
If $\dim (U_{X- x}) +1\geq \dim (U)$ for any $x\in X$, then by induction
$$\dim (U)\leq \dim (U_{X-x}) +1\leq |X-x| -2 +1=|X|-2.$$
In this case, we are done.

In the other case there exists an $x\in X$ such that  $\dim (U_{X-x}) +2\leq \dim (U)$. Suppose $C$ is a cyclic ordering such that $u^C\in U$. Then we may write $C=[x_1,\ldots, x_n]$ with $x_n=x$. We put $\infty=x_1$ and fix the linear ordering $\infty=x_1<x_2<\cdots<x_n=x$ of $X$.
Let  $$T:=\{(\infty, y,z)\mid y,z\in X, y<z\}, T':= \{(\infty, y,z)\mid y,z\in X-x, y<z\}, $$
Then we have $0=u^{C'}_{|T}\in U_{|T}$, where $C'$ arises by reversing $C$.  This makes $U_{|T}$ a linear subspace, and since $\dim(U_{|T'})+2=\dim (U_{X-x}) +2\leq \dim (U)=\dim(U_{|T})$, there must exist two vectors $v,w$ parallel to $U$ such that $v_{|T'}=w_{|T'}=0$, and such that $v_{|T}, w_{|T}$ are linearly independent. Then there must exist $(\infty, a,x), (\infty,b,x)\in T\setminus T'$ so that the restrictions of $v,w$ to $\{(\infty, a,x), (\infty,b,x)\}$ are already linearly independent, and we may assume $v(\infty,a,x)=1\neq w(\infty,a,x)$, $v(\infty,b,x)\neq 1= w(\infty,b,x)$, and $a<b$. Then $v_{\{\infty,a,b,x\}|T}$ is not cyclic, and hence $v_{|T}$ is not cyclic, contradicting that $v_{|T}=v_{|T}+u^{C'}_{|T}\in U_{|T}$.

So $\dim(U)\leq |X|-2$. We now prove the second claim, that $\dim(U)= |X|-2$ if and only if $U=U(G)$ for some binary tree on $X$. Sufficiency being easy, we only prove necessity. For that, we apply  induction on $|X|$ again. For the induction step, we need to show that if $U_{X-x}=U(G')$ for a binary tree $G'$ on $X-x$ and $\dim(U)=\dim(U_{X-x})+1$, then there is a binary tree $G$ so that $U=U(G)$.

We consider a cyclic ordering $C$ so that $u^C\in U$ again, and the bases $T,T'$ of $\UU^X$ resp. $\UU^{X-x}$ as before. This time  $\dim(U_{|T'})+1=\dim (U_{X-x}) +1\leq \dim (U)=\dim(U_{|T})$, so that there is a unique nonzero vector $v$ parallel to $U$ so that $v_{|T'}=0$. Let $Y:=\{y\in X\mid v(\infty,x,y)=1\}$. If $u^{C_1}\in U$ for a cyclic ordering $C_1$, then $u^{C_1}+v=u^{C_2}$ for some other ordering $C_2$. Then
$$\{C_1,C_2\}=\{[\infty, x_1,\dots, x_k, y_1,\ldots, y_l, x], [\infty, x_1,\dots, x_k, x, y_1,\ldots, y_l]\},$$
where $Y=\{y_1,\ldots, y_l\}$. It follows that $Y$ is consecutive in each cyclic ordering $C_1$ such that $u^{C_1}\in U$, and hence that $Y$ is consecutive in each cyclic ordering consistent with $G'$, and as $G'$ is a tree, this implies that $Y$ is a split of $G'$, say $G'-e$ has a component $W$ with $W\cap (X-x)=Y$. Then replacing $e$ with two series edges with common point $v$ and adding an edge $xv$, we obtain a graph $G$ so that $U=U(G)$.
\endproof

\section{\label{sect:algorithm} An algorithm for determining a level-1 network from a set of quartets}
If $Q$ is a set of quartets on a set of taxa $X$, the we say that a level-1 network $G$ {\em displays} $Q$ if  $Q\subseteq\QQ(G)$, and that $Q$ is {\em inconsistent} if there is no level-1 network $G$ that displays $Q$.
By Lemma \ref{lemma:char2} and Theorem \ref{thm:char2}, a level-1 network $G$ displays a quartet set $Q$ if and only if $U(G)\subseteq U(Q)$.

By Theorem \ref{thm:dim}, if $U(Q)$ is cyclic, 
then $\dim U(Q)_Z\leq |Z|-2$ for all $Z\subseteq X$. We say that $Q$ is {\em insufficient} if $U(Q)$ is not cyclic, so that $\dim U(Q)_Z>2$ for some $Z\in \binom{X}{4}$, and a we call such a $Z$ a {\em witness} for the insufficiency of $Q$.

In this section, we describe an algorithm that solves the following problem in $O(n^6)$ time.
\begin{tabbing}
{\bf Given:} \= A set of quartets $Q$ on a set of $n$ taxa $X$, and a taxon $\infty\in X$.\\
{\bf Find:}\> A binary level-1 network $G$ on $X$ that displays $Q$, such that $|\Sigma(G)|$ is as large as possible, or \\
		\>an inconsistent subset $Q'\subseteq Q$ with $|Q'|\leq \binom{n-1}{2}$, or \\
		\>a witness $Z\in \binom{X}{4}$ such that $\infty\in Z$.
\end{tabbing}
If the input set of quartets $Q$ is 
dense, 
then $U(Q)$ is cyclic by Lemma \ref{theorem:dense} and hence $Q$ is sufficient. So for dense $Q$, the algorithm presented here finds a binary level-1 network $G$ on $X$ that displays $Q$ or it outputs a short certificate of inconsistency of $Q$. At the end of this section, we will also consider the special case where $Q=\{q\in \QQ(G)\mid\infty\in \underline{q}\}$. We will argue that then, a variant of our algorithm reconstructs $G$ from $Q$ in $O(n^3)$ time.

Without any assumptions on the input set $Q$, our algorithm is heuristic. If $Q$ is ambiguous, then the algorithm may fail to produce a network or to detect inconsistency of the input, but outputs a witness $Z$ instead. But any witness $Z$ also points out how the insufficiency of $Q$ may be repaired. For any quartet $q$ such that $\underline{q}=Z$, we have $\dim U(Q+q)_Z=2<3=\dim U(Q)_Z$, and hence $\dim U(Q+q)< \dim U(Q)$. This suggests an interactive use of our algorithm in a practical setting: if a witness $Z$ is detected, find some new quartet $q$ such that $\underline{q}=Z$ and add it to $Q$, and repeat this until until the algorithm terminates with a $G$ displaying $Q$ or an inconsistent $Q'$. As $\dim U(Q)\leq \dim \UU^X\leq O(n^2)$ and $\dim U(Q)$ decreases with each addition to $Q$, at most $O(n^2)$ new quartets will be needed to force one of these outcomes.

\subsection{\label{ss:la}Constructing $U(Q)$}
We are given a set of quartets $Q$ on $X$, and a taxon $\infty\in X$. We fix a linear ordering $<$ of $X$, and we consider $T:=\{(\infty, a,b) \mid a<b\}$. Recall that
$$U(Q):=\{u\in\UU\mid u(a,b,c)+u(a,b,d)=0\text{ for all }(ab|cd)\in Q\}.$$
Hence
$$U(Q)_{|T}=\{x\in \GF(2)^T\mid a^qx=b^q\text{ for all }q\in Q\},$$
where $a^qx=b^q$ is a sparse linear equation, as described in subsection \ref{ss:basis}.

We may determine $u, V$ such that
$U(Q)_{|T}=\{u+Vy\mid y\in \GF(2)^d\}$ by applying 
the algorithm described in subsection \ref{ss:linalg}.
The number of variables is then $m:=|T|=\binom{n-1}{2}=O(n^2)$ and the number of equations is $k:=|Q|=O(n^4)$.
The resulting running time is $O(km+m^3)=O(n^6)$, as required.
At termination, the set $K$ determined by the linear algebra procedure corresponds to a set $Q'\subseteq Q$ with $|Q'|=O(n^2)$ such that $U(Q)=U(Q')$. If it turns out that $U(Q)=\emptyset$, then $U(Q')=\emptyset$. In that latter case, we may exit quoting $Q'$ as an inconsistent subset of $O(n^2)$ quartets.

\subsection{\label{ss:badZ}Finding the subspace $U$ of all cyclic vectors in $U(Q)$, or a witness $Z$}
Let
\begin{equation}\label{Ucycl} U:=\{u\in U(Q)\mid u\text{ is cyclic }\}.\end{equation}
We describe how to find, either
\begin{enumerate}
\item a set of quartets $Q^+$ so that $U(Q'\cup Q^+)=U$, or
\item a witness $Z\in \binom{X}{4}$ such that $\infty\in Z$.
\end{enumerate}
in $O(n^5)$ time.

We initialize the algorithm with $Q^+=\emptyset$ and consider the elements of $\mathcal{Z}_\infty:=\{Z\in \binom{X}{4}\mid \infty\in Z\}$ one by one. If $U(Q')_Z$ contains only cyclic vectors, then we do nothing. If $U(Q')_Z$ contains a non-cyclic vector $w$, then there are three cases depending on the dimension of $U(Q')_Z$:
\begin{enumerate}
\item $\dim (U(Q')_Z)=1$ and $U(Q')_Z=\{w, w+{\bf 1}_Z\}$; or
\item $\dim (U(Q')_Z)=2$ and $U(Q')_Z=\{u^C, u^C+{\bf 1}_Z, w, w+{\bf 1}_Z\}$, where $C$ is a cyclic ordering of $Z$; or
\item $\dim (U(Q')_Z)=3$.
\end{enumerate}
In case (1), there are no cyclic vectors in $U(Q')$, and hence there are no cyclic vectors in $U(Q')$. Then we terminate quoting $Q'$ as an inconsistent set of quartets.
In case (2), we may construct a quartet $q$ such that $\underline{q}=Z$ and such that  $U(Q'+q)_Z=\{u^C, u^C+{\bf 1}_Z\}$, e.g. if $C=[a,b,c,d]$, then taking $q=(ab|cd)$ will do. Note that this implies that $U(Q'+q)_Z\supseteq U$. We then add this quartet to $Q^+$.
In case (3), $Z$ is a witness of the insufficiency of $Q'$, and hence a witness of the insufficiency of $Q$. We then terminate the algorithm outputting $Z$.

After dealing with each $Z\in \mathcal{Z}_\infty$, we have obtained a $Q^+$ such that $U\subseteq U(Q'\cup Q^+)$ and such that $U(Q'\cup Q^+)_Z$ contains only cyclic vectors for each $Z\in \mathcal{Z}_\infty$. By Lemma \ref{lemma:noncyclic}, all vectors of $U(Q'\cup Q^+)$ are cyclic, and so $U=U(Q'\cup Q^+)$.

Given a fixed $Z\in \binom{X}{4}$, it takes $O(\dim(U(Q'))$ time to find the vector representation of $U(Q'')_Z$
by restricting the vector $u'$ and the columns of $V'$ to $T\cap(\infty\times Z\times Z)$. Note that each restricted vector has constant length, so that it takes $O(\dim(U))$ time to remove duplicate vectors.
Then, it takes $O(1)$ time to compute $\dim(U(Q')_Z)$ and to determine whether $U(Q')_Z$ contains non-cyclic vectors.
Since we have to consider $O(n^3)$ sets $Z$, the overall time is $O(n^5)$.

In what follows, we assume that we have $u,V$, and $Q''\subseteq Q'\cup Q^+$ such that
$$U=U(Q'')=U(Q'\cup Q^+)=\{u+Vy\mid y\in GF(2)^d\},$$
and $|Q''|+d=|T|$. It takes $O(n^6)$ time to obtain such $u,V, Q''$ from $Q'\cup Q^+$, as in subsection \ref{ss:la}.

\subsection{Finding a cyclic ordering $C$}\label{ss:cyclic}
Given any $u\in U$, it takes $O(n^2)$ time to construct a cyclic ordering $C$ so that $u=u^C$. Just proceed by adding the elements of $X$ one by one, starting with $\infty$ and any two other elements. If a set $Y=\{\infty, y_1,\ldots, y_s\}$ has been cyclically ordered as $C=[\infty, y_1,\ldots, y_s]$ so that
$u^C=u_Y$, consider an $x\in X\setminus Y$. If the set $Y':=\{y\in Y-\infty\mid u(\infty, y,x)=1\}$ is such that $Y'=\{y_1,\ldots, y_t\}$ for some
$t$, then we may extend $C$ with $x$ to $C'=[\infty, y_1,\ldots, y_t, x, y_{t+1},\ldots, y_s]$, so that  $u^{C'}=u_{Y+x}$. If not, then there are $i,j<s$ so that $i<j$ and
$u(\infty, y_i, x)=0, u(\infty, y_j, x)=1$. Then with $Z=\{\infty, y+i, y_j, x\}$, the restriction $u_Z$ is not cyclic. This would contradict the choice of $u$, as each $u\in U$ is cyclic.

\subsection{Finding a set of splits $\Sigma$} Let $C$ be the cyclic ordering obtained in the above, and let $U$ be the affine set of cyclic vectors from $U(Q)$. We want to construct the set $\Sigma(U)$.
We claim that each $S\in \Sigma(U)$ is consecutive in $C$. If not, then there are $a,b,c,d\in X$ so that $C\succeq [a,b,c,d]$ and $a,c\in S$, $b,d\not\in S$. Then with $Z=\{a,b,c,d\}$, the vector $(u^C+v^S)_Z\in U$ is not cyclic. This contradicts that each $u\in U$ is cyclic.
Thus we have
$$\Sigma(U)=\{S\subseteq X\mid |S|, |X\setminus S|\geq 2, ~S\text{ is consecutive in }C, ~ v^S\text{ parallel to }U\}.$$
This makes $\Sigma(U)$ easier to compute. There is a straightforward way taking $O(n^4)$ time, as follows.
Enumerate the $O(n^2)$ sets $\{S\subseteq X\mid |S|, |X\setminus S|\geq 2, ~S\text{ is consecutive in }C\}$ in as much time.
For each fixed $S$ it takes $O(n^2)$ time to decide if $v^S$ is parallel to $U$, since this is the case if and only if
$a^qv^S=0$ for all $q\in Q''$, each $a^q$ is sparse and $|Q''|$ is $O(n^2)$

The set $\Sigma(U)$ is {\em cross-free}, that is, for each $S,T\in \Sigma(U)$ we have
$$S\subseteq T, T\subseteq S, S\cup T=X, \text{ or }S\cap T=\emptyset.$$
 If not, there would be  $S,T\in \Sigma(U)$ and $a,b,c,d$ such that $C\succeq [a,b,c,d]$ and $a\in S\cap T$, $b\in S\setminus T$, $c\not\in S\cup T$, $d\in T\setminus S$. Then with $$U':=\{u^C+\lambda_1 {\bf 1}_X+ \lambda_2 v^S+\lambda_3 v^T\mid \lambda_i\in \GF(2)\},$$
we have $U'\subseteq U$ and $\dim U'_Z=3$, where $Z=\{a,b,c,d\}$. This contradicts that $U$ contains only cyclic vectors.

Since $u^C\in U$ and $v^S$ is parallel to $U$ for each $S\in \Sigma(U)$, the set
$$U_0:=\{u^C+\sum_{S\in \Sigma(U)} \lambda_S v^S\mid \lambda_S\in \GF(2)\}$$
is contained in $U$.
\subsection{Constructing $G$}\label{ss:constrG}
Given $C$, it takes $O(n)$ time to build an array that gives the index of each taxon on $C$ in clockwise order.
 Then obtaining the index of any taxon in the cyclic order $C$ takes constant time.
 We may also assume that each set $S\in \Sigma(U)$ is given to us as an interval in $C$, specified by its first and last index in clockwise order. Then testing containment of two such intervals takes $O(1)$ time.
  The tree-representation of the cross-free set
$\Sigma(U)$ of cardinality $O(n)$ can then conveniently be obtained in $O(n^2)$ time from a circular arrangement of the taxa in the order
suggested by $C$,
and when each vertex of degree $d>3$ in this tree is replaced by a circuit of length $d$ such that the ordering of the (groups of) taxa around this circuit respects $C$, a level-1 network $G$ is obtained such that
$\Sigma(G)=\Sigma(U)$, and such that $C\in \CC(G)$. It follows that $U_0=U(G)$. The construction of this $G$ takes $O(n^2)$ time.

Suppose that $G'$ is such that $U(G')\subseteq U(Q)$, i.e. $G'$ displays the quartets  from $Q$. Then $U(G')\subseteq U$. We now have $\Sigma(G')\subseteq \Sigma(U)$, since if $S\in \Sigma(G')$, then $v^S$ is parallel to $U(G')$, and hence $v^S$ is parallel to $U$, so that $S\in \Sigma(U)$.
Therefore, $G$ as constructed above with $\Sigma(G)=\Sigma(U)$ has the maximum number of splits among all level-1 networks displaying the quartets from $Q$, or briefly
$$|\Sigma(G)|=\max\{|\Sigma(G')| \mid U(G')\subseteq U(Q)\},$$
as required.

\subsection*{Remark} It is possible that $\dim(U_0)<\dim(U)$ in the final stage of the algorithm. Upon decomposing along the splits $\Sigma(U)$, one will then find that at least one of the remaining split-free components has dimension $>1$.

\subsection{A refined algorithm}
Let us denote the set of quartets of a level-1 network $G$ that contain a fixed taxon $\infty \in X$ by $Q_{\infty}(G)$. So
\[Q_{\infty}(G):=\{q\in \QQ(G)\mid \infty \in \underline{q}\}\]
Now suppose that we are given a set of quartets $Q$, such that each quartet from $Q$ contains  $\infty$.
Such a set of quartets is called \begin{em} level-1 like \end{em} if $Q=Q_{\infty}(G)$ for some level-1 network $G$.
In this section we will argue that there is a $O(n^3)$ 
algorithm to solve the following problem.

\begin{tabbing}
{\bf Given:} \= A set of quartets $Q$ with $\infty \in \underline{q}$ for each $q\in Q$ that is level-1 like 
\\
{\bf Find:}\> A level-1 network $G$ such that $Q=Q_{\infty}(G)$.
\end{tabbing}

Assume that $Q=Q_infty(G)$. Then
 $U(Q)=U(G)$, by Theorem \ref{thm:rooted}, and  $G$ is 
unique 
by Theorem
\ref{thm:unique}. Note furthermore that it suffices to output a $G'$ that displays $Q$ and has $|\Sigma(G')|$ as large as possible,
since in that case $U(G')\subseteq U(Q)=U(G)$ and by Lemma \ref{lemma:dimG}, $\dim (U(G')) =|\Sigma(G')|/2+1\geq |\Sigma(G)|/2+1=\dim(U(G))$. It follows that $U(G')=U(G)$, and hence $G'=G$.

Let $ <$ be a linear ordering of $X$. 
We consider a signed auxiliary graph $H=H(Q,\infty)$ with
vertex set 
 	\[T:=\{ (\infty, a,b) \mid a, b \in X\setminus \{\infty\}, a<b \}, \mbox{ and }\]
 and edge set $E=E_0\cup E_1$ consisting of even edges
 \[ E_0:=\{ \{(\infty, a,b),(\infty, a,c)\} \mid (\infty a | b c ) \in Q, \; a<b<c \}\cup \{\{(\infty,a,c),(\infty,b,c)\}:(\infty c|ab)\in Q, \; a<b<c\}\] and odd edges
 \[ E_1:=\{ \{(\infty, a,b),(\infty, b,c)\} \mid (\infty b | a c ) \in Q, \; a<b<c \}.\]
Let \[ \Sigma:=\{ S(W)  :  W \mbox{ the vertex set of a component of } H \},\] where
 \[S(W)= \bigcup \{ \{a,b\} \mid (\infty, a,b) \in W \}.\]
From the signed graph $H$, any $u\in U(Q)=\{ u\in \UU : u(\infty, a,b)+u(\infty, a,c)=0 \mbox{ for all }(\infty a | b c ) \in Q \}$ can be constructed by choosing for every vertex set of a component $W$ of $H$ a fixed value $u_t\in \{0,1\}$ for some arbitrary $t\in W$, and then putting $u_s=u_t$ for every $s$ in that same component $W$ at even distance from $t$, and putting $u_s\neq u_t$ for every $s$ at odd distance from  $t$. 
In other words, for any fixed $u\in U(Q)$:
\begin{equation}\label{V'}
U(Q)
 = \{u+ V' \alpha \mid  \alpha \in GF(2)^d \}
 \end{equation}
  where $V'$ is a matrix with a column 
  $v'$ parallel to $\UU$ with 
  $v'_{|T} = \chi^{W}$,  for every vertex set $W$ of a component 
  of  $H$.

\begin{lemma}\label{lemma:quick} $\Sigma
=\{S\in \Sigma(G): \infty \notin S\}\cup\{X-\infty\}$.
\end{lemma}
\proof
By Lemma \ref{lemma:dimG},
 	\[U(Q)=U(G)=\{u+ V\alpha| \alpha \in GF(2)^d\},\]
where 
$V$ is a matrix with a column 
$v^S$ for every nontrivial 
split $S$ of  $G$ such that $\infty \notin S$, and an all-one column $1=v^{X-\infty}$.
In other words, $V$ has a column $v^S$ for each $S\in {\mathcal S}$, where ${\mathcal S}$ is the following laminar collection of splits:
\[{\mathcal S}=\{S\in \Sigma(G)\mid \infty \notin S\}\cup\{X-\infty\}.\]
By (\ref{V'}),
every column of $V$ is a linear combination of columns of $V'$.
Because for every column $v'$ of  $ V'$ it holds that $v'_{|T} = \chi^{W}$, i.e. $\supp(v')\cap T= W$ for some vertex set $W$ of a component of $H$, 
 the supports of columns 
 $v'$ of  $V'$ are disjoint on $T$. Thus, to realize a column 
 $v^S$ of  $V $ as 
 a linear combination of such $v'$, we need to take 
\[ v^S=\sum  \{ v' \mbox{ column of } V' \mid v'_{|T}=\chi^{W}  \mbox{ for a component $W$ contained in } \supp(v^S) \}.\]

We claim that for every  
$S\in {\mathcal S}$, there is a
component $W$ of $H$ such that $S=S(W)$.
To see this, let $S\in {\mathcal S}$, and consider the inclusionwise maximal splits
$S_1,\ldots, S_k$  that are contained in  
$S$ (possibly, some of these 
$S_i$ are trivial splits, otherwise they are again contained in ${\mathcal S}$). 
Note that $S$ is the disjoint union of these 
$S_i$, and that $k>1$. 
We have
\[v^S
= \sum \{ v' \mid v'{_T}=\chi^{W} \mbox{ for a $ W$ contained in $\supp(v^S)$ but not in
 any } \supp(v^{S_i}) \}
                      + v^{S_1} + .. + v^{S_k}.\]

Since $k>1$, 
there is at least one component $W$ such that $W$ is contained in
$\supp(v^S)$, but not in any 
$\supp (v^{S_i})$. Let $W_S$ denote such a component.
  Then  the map
 $S\mapsto W_S$ from ${\mathcal S}$ to the set of components of $H$ is injective because ${\mathcal S}$ is laminar.
 There are $d$ components in $H$, where $d= \dim(U(Q))$, and there are 
 also $d$  splits $S$ in the collection ${\mathcal S}$.
  Therefore the map $S\mapsto W_S$ is in fact a bijection, and it follows that for each $S\in {\mathcal S}$ there is a unique $W_S$ such that
  $W_S$ contained in $\supp(v^S)$ but not in any $\supp(v^{S_i})$.

 Now we prove that $S=S(W_S)$ for each $S\in {\mathcal S}$. Let $a\in S$. Since $k>1$ we can choose $b\in S$
 such that $a$ and $b$ are in  distinct $S_i$. Then 
 $(\infty,a,b)\in W_S$ 
 since $(\infty,a,b)$ is in $\supp(v^S)$ and not in any 
 $\supp(v^{S_i})$, so $a\in S(W_S)$.
 Conversely, let $a\in S(W_S)$, then $(\infty,a,b)\in W_S$ for some $b>a$ or $(\infty, c,a )\in W_S$ for some $c<a$. But then
  $(\infty,a,b)$ (or $(\infty,c,a)$) is in $\supp(v^S)$, so $a,b\in S$ (or $a,c\in S$). In both cases $a\in S$.

Since $|\Sigma|=d=|{\mathcal S}|$, and by the above claim ${\mathcal S}\subseteq
\Sigma$, we conclude
$\Sigma={\mathcal S}$.
\endproof

We derive an algorithm with running time $O(n^3)$ for the problem posed at the beginning of this subsection.
Reading the data and constructing  the graph $H$ takes  $O(|T|+|Q|)= O(n^3)$ time.
Finding the components of $H$ takes time   $O(|T|+|Q|)= O(n^3)$.
Now some fixed $u\in U(Q)$ can be constructed in time $O(|T|)=O(n^2)$, and a cyclic order $C$ such that $u=u^C$ can be constructed as in
Subsection \ref{ss:cyclic} in time $O(n^2)$. We may then assume to have a data structure in which the rank number of any taxon in this cyclic order can be determined in constant time.
 Determining $S(W)$ (finding the first and the last rank number in the cyclic order $C$ of taxa in  $S(W)$)then takes time 
     $O(|V(W)|)$ for each 
     $W$,
    so 
     determining $\Sigma$, i.e. the set of splits of $G$, takes time 
     $O(|T|)=O(n^2)$.
      Constructing a network $G$ from its set of splits takes time $O(n^2)$ as in Subsection \ref{ss:constrG}. 
      So the entire algorithm takes $O(n^3)$ time in this special case.

      Since the bottleneck for the complexity  of this special case seems to be the reading of the quartet data, we conjecture that a faster algorithm would be possible if the quartets were given by an oracle.

\section{\label{sect:outtro} Discussion }
\subsection{Summary}
We described an algorithm which, given a set of quartets $Q$ on a set $X$ of $n$ taxa, either finds a compatible level-1 network $G$ with a maximum number of splits, or exits with information that can be used to augment $Q$, preventing further failures. It has been noted 
that if we input a dense quartet set $Q$, i.e. a $Q$ such that $\{\underline{q}\mid q\in Q\}=\binom{X}{4}$, then the process 
always terminates with either a level-1 network or an inconsistent subset of $Q$ of cardinality $O(n^2)$.
This settles an open problem posed by Gambette, Berry and Paul \cite{Gambette2012}. Their question whether it
is possible to reconstruct in polynomial time a (simple) unrooted level-1 network from a dense quartet set, is answered affirmatively.

\subsection{Implementation}
The running time of $O(n^6)$ of our general algorithm clearly makes that there will be a quite sharp upper bound on the number of taxa beyond which this method will no longer be practically feasible. We estimate that our method will feasible for problems involving hundreds of taxa. The bottleneck of our algorithm is the linear algebra over $\GF(2)$. Fortunately, present-day computers are simply made for doing linear algebra over $\GF(2)$, as vectors in $\GF(2)^{64}$ may naturally be stored a machine word, and addition of such vectors is an elementary operation for the CPU (XOR). Therefore, the hidden constant in the stated $O(n^6)$ running time will be very modest when applying the linear algebra as described. The problem of solving systems of linear equations over $\GF(2)$ also arises in the factorization of integers, and perhaps because of this it has attracted the attention of many algorithm designers \cite{Wiedemann1986, Coppersmith1994}. In particular, Wiedemann \cite{Wiedemann1986} considers the problem of solving sparse linear equations, and we note that the system of equations $Ax=b$ occurring in the method presented here is clearly sparse with at most 4 nonzero entries per row of $A$. These, and other methods have been implemented in \cite{m4ri}.

At any time during the execution of our algorithm for solving $k$ equations in $m$ variables over $\GF(2)$, the storage of at most $m$ vectors over $\GF(2)$ of length $m$ takes at most $m^2$ bits. With 2 Gigabites of RAM ($=2^{31}$ bytes $= 2^{34} $ bits), a value of $m=2^{17}=131072$ is feasible. In our application to reconstructing a phylogenetic network on a set of $n$ taxa, the value of $m$ will be  about $n^2/2$. This implies that up to some $500$ taxa, the memory requirements of the algorithm should not be a bottleneck.

An important part of the algorithm described in this paper has been implemented by Willem Sonke  \cite{WillemSonkeImpl}, with promising results. The implementation and some computational results are described in his bachelor thesis
\cite{WillemSonkeThesis}.

\subsection{Quartet inference}

 If $G$ is a tree on $X$, then the following {\em dyadic inference rules} apply:
 \begin{enumerate}
\item if $(ab|cd), (ab|de)\in\QQ(G)$, then $(ab|ce)\in \QQ(G)$, and
\item if $(ab|cd), (ac|de)\in \QQ(G)$, then $(ab|ce)\in \QQ(G)$.
\end{enumerate}
An $O(n^5)$ algorithm for finding the quartets that follow from a set of quartets $Q$ by applying dyadic inference to exhaustion is described in \cite{ErdosSteelSzekelyWarnow1999}.

The situation seems to be less straightforward for level-1 networks. A finite case analysis reveals that if $G$ is a level-1 network, then
\begin{enumerate}
\item if $(ab|cd), (ab|de)\in\QQ(G)$, then $(ab|ce)\in \QQ(G)$, and
\item if $(ab|cd), (ac|de)\in \QQ(G)$, then $\QQ(H_1)\subseteq \QQ(G)$ or $\QQ(H_2)\subseteq \QQ(G)$.
\end{enumerate}
Here $H_1, H_2$ are the level-1 networks in Figure \ref{fig:dyadic}. Note that  $(ab|ce)\not\in \QQ(H_2)$, so that the second dyadic inference rule for trees does not apply to level-1 networks.
Due to the disjunction implied in the second rule, it is not obvious to us how to even define the closure of a given set of quartets under these inference rules, let alone compute it in polynomial time.

The linear algebra methods of the present paper inspire the notion of quartet inference where a quartet $q$  is considered to be implied by a set of quartets $Q$ if $U(Q+q)=U(Q)$.
\section{Acknowledgement}
The authors would like to thank Steven Kelk
for initiating this research by posing the question of Gambette, Berry and Paul \cite{Gambette2012} to us.

\bibliographystyle{plain}
\bibliography{quartets}

\end{document}